\documentclass[11pt]{amsart}

\def\CC{{\mathcal C}}

\def\EE{{\mathcal E}}

\def\FF{{\mathcal F}}

\def\PP{{\mathcal P}}
\def\PPp{{\mathcal P'}}

\newcommand{\B}{{\ensuremath{\mathbb  B}}}

\newcommand{\BBn}{{\ensuremath{{\B}^{n+1}}}}
\newcommand{\BBnbar}{{\ensuremath{{\overline{\B}^{n+1}}}}}
\newcommand{\BBtre}{{\ensuremath{{{\B}^{3}}}}}

\newcommand{\C}{{\ensuremath{\mathbb  C}}}
\newcommand{\Cbar}{\ensuremath{\overline{\mathbb  C}}}

\newcommand{\CPone}{\ensuremath{{{\mathbb   {CP}}^1}}}
\renewcommand{\d}{\ensuremath{\operatorname{d}}}

\newcommand{\deta}{\ensuremath{\operatorname{d\eta}}}

\newcommand{\dmu}{\ensuremath{\operatorname{d\mu}}}

\newcommand{\dxv}{\ensuremath{d\mathbf{x}}}

\newcommand{\dzeta}{\ensuremath{\operatorname{d\zeta}}}
\newcommand{\D}{\ensuremath{{\mathbb   D}}}

\newcommand{\e}{\ensuremath{{\operatorname{e}}}}
\newcommand{\eps}{{\epsilon}}
\newcommand{\epsv}{\ensuremath{{\underline{\epsilon}}}}
\newcommand{\ev}{\ensuremath{{\mathbf{e}}}}

\newcommand{\End}{\ensuremath{\operatorname{End}}}

\newcommand{\isom}{\ensuremath{\operatorname{\approx}}}

\newcommand{\J}{\ensuremath{{\operatorname{J}}}}
\newcommand{\Jac}{\ensuremath{{\operatorname{Jac}}}}

\newcommand{\mapfromto}[3]{\hbox{\ensuremath{#1 : #2 \longrightarrow #3}}}

\newcommand{\ooo}{\ensuremath{{\emph{o}}}}

\newcommand{\rv}{\ensuremath{{\mathbf{r}}}}
\newcommand{\R}{\ensuremath{\mathbb  R}}

\newcommand{\Rplus}{\ensuremath{{{\mathbb  R}_+}}}

\newcommand{\Sm}{\ensuremath{{\setminus}}}
\newcommand{\Sen}{\ensuremath{{{\mathbb  S}^1}}}

\newcommand{\Sn}{\ensuremath{{{\mathbb  S}^n}}}
\newcommand{\Snminus}{\ensuremath{{{\mathbb  S}^n_-}}}
\newcommand{\Snplus}{\ensuremath{{{\mathbb  S}^n_+}}}

\newcommand{\Sto}{\ensuremath{{{\mathbb  S}^2}}}
\newcommand{\Stre}{\ensuremath{{{\mathbb  S}^3}}}

\newcommand{\vv}{\ensuremath{{\mathbf{v}}}}

\newcommand{\Vol}{\ensuremath{{\operatorname{Vol}}}}

\newcommand{\wf}{\ensuremath{\widehat f}}

\newcommand{\wv}{\ensuremath{{\mathbf{w}}}}

\newcommand{\X}{\ensuremath{\mathbb  X}}
\newcommand{\xv}{\ensuremath{{\mathbf{x}}}}
\newcommand{\Y}{\ensuremath{\mathbb  Y}}
\newcommand{\zbar}{\ensuremath{{\overline{z}}}}
\newcommand{\zerov}{\ensuremath{{\mathbf{0}}}}
\newcommand{\zetav}{\ensuremath{{\underline{\zeta}}}}
\newcommand{\zv}{\ensuremath{{\mathbf{z}}}}

\def\eps{\epsilon}

\def\R{\mbox{$\mathbb R$}}

\def\C{\mbox{$\mathbb C$}}

\def\D{\mbox{$\mathbb D$}}

\def\Mobius{M{\"o}bius}
\def\Poincare{Poincar{\'e}}
\newtheorem{newthm}{Theorem}

\newtheorem{theorem}{Theorem}
\newtheorem{lemma}[theorem]{Lemma}
\newtheorem{proposition}[theorem]{Proposition}
\newtheorem{corollary}[theorem]{Corollary}
\newtheorem{conjecture}{Conjecture}
\newtheorem{remark}[theorem]{Remark}
\newtheorem{defthm}[theorem]{Definition and Theorem}
\newtheorem{definition}[theorem]{Definition}
\newcommand{\ALIGN}{\begin{align*}}
\newcommand{\ENDALIGN}{\end{align*}}
\newcommand{\ENUM}{\begin{enumerate}}
\newcommand{\ENUMa}{\begin{enumerate}[a.]}
\newcommand{\ENUMi}{\begin{enumerate}[i)]}
\newcommand{\ENDENUM}{\end{enumerate}}
\newcommand{\ITMZ}{\begin{itemize}}
\newcommand{\ENDITMZ}{\end{itemize}}
\newcommand{\EQN}[1] { \begin{equation}\label{#1} }
\newcommand{\ENDEQN}{\end{equation}}
\newcommand{\THM}{\begin{theorem}}
\newcommand{\REFEXA}[1] { \begin{example}\label{#1} }
\newcommand{\ENDEXA}{\end{example}}
\newcommand{\REM}{ \begin{remark}}
\newcommand{\ENDREM}{\end{remark}}
\newcommand{\REFTHM}[1] { \begin{theorem}\label{#1} }
\newcommand{\RREFTHM}[2] { \begin{theorem}[#1]\label{#2} }
\newcommand{\ENDTHM}{\end{theorem}}
\newcommand{\REFNTH}[1] { \begin{newthm}\label{#1} }
\newcommand{\ENDNTH}{\end{newthm}}
\newcommand{\REFPROP}[1]{\begin{proposition}\label{#1} }
\newcommand{\RREFPROP}[2]{\begin{proposition}[#1]\label{#2} }
\newcommand{\PROP}{\begin{proposition}}
\newcommand{\ENDPROP}{\end{proposition} }
\newcommand{\REFDEF}[1]{\begin{definition}\label{#1} }
\newcommand{\DEF}{\begin{definition}}
\newcommand{\ENDDEF}{\end{definition} }
\newcommand{\REFLEM}[1]{\begin{lemma}\label{#1} }
\newcommand{\RREFLEM}[2]{\begin{lemma}[#1]\label{#2} }
\newcommand{\LEM}{\begin{lemma}}
\newcommand{\ENDLEM}{\end{lemma} }
\newcommand{\REFCOR}[1]{\begin{corollary}\label{#1} }
\newcommand{\COR}{\begin{corollary}}
\newcommand{\ENDCOR}{\end{corollary}}
\newcommand{\REFCONJ}[1]{\begin{conjecture}\label{#1} }
\newcommand{\CONJ}{\begin{conjecture}}
\newcommand{\ENDCONJ}{\end{conjecture}}
\newcommand{\RMRK}{\begin{remark}}
\newcommand{\ENDRMRK}{\end{remark}}
\newcommand{\REFDEFTHM}[1] { \begin{defthm}\label{#1} }
\newcommand{\ENDDEFTHM}{\end{defthm}}

\newcommand{\lemref}[1]{Lemma~\ref{#1}}

\newcommand{\propref}[1]{Proposition~\ref{#1}}

\newcommand{\PROOF}{\begin{proof}}
\newcommand{\ENDPROOF}{\end{proof}}

\usepackage{geometry}                
\usepackage{graphicx}
\usepackage{amssymb}
\usepackage{epstopdf}
\usepackage[active]{srcltx}

\title{Conformally Natural extensions revisited}
\author{Carsten Lunde Petersen}
\date{\today}                                           

\begin{document}
\maketitle
\begin{abstract}
In this note we revisit the notion of conformal barycenter of a measure on $\Sn$ as 
defined by Douady and Earle \cite{D-E}. The aim is to extend rational maps 
from the Riemann sphere $\Cbar\isom\Sto$ to the (hyperbolic) three ball $\BBtre$ 
and thus to $\Stre$ by reflection. 
The construction which was pioneered by Douady and Earle in the case of homeomorphisms 
actually gives extensions for more general maps such as entire transcendental maps 
on $\C\subset\Cbar$. 
And it works in any dimension. 
\end{abstract}
\section{Intoduction}
Let $G = G_n$ denote the group of {\Mobius} transformations of {\Mobius} space
$\widehat\R^n=\R^n\cup\{\infty\}$ preserving the $n$-sphere $\Sn$: 
$$
\Sn = \{(x_1, \ldots x_{n+1})\in \R^{n+1} | x_1^2 + \ldots x_{n+1}^2 = 1 \}.
$$
as well as the enclosed ball $\BBn$. 
Then each element $g\in G$ acts on $\BBn$ as a hyperbolic isometry, 
that is it preserves the Riemannian metric $2|\dxv|/(1-|\xv|^2)$.
Moreover each $g$ is conformal and thus also 
acts as a conformal automorphism of both $\Sn$ and $\BBn$. 
Mostow \cite{Mostow} proved that any conformal isomorphism of $\BBn$ 
and/or $\Sn$ is an element of $G$, so that we may also define $G$ as the 
conformal automorphism group of $\BBn$ and/or $\Sn$.
We let $G_+=G_{n,+}$ denote the index two subgroup consisting of orientation preserving 
conformal automorphisms. 
And we let $c$ denote the reflection in the coordinate plane $x_{n+1} =0$,  
so that $G$ is generated by $G_+$ and $c$, i.e.~$G = <G_+, c>$. 

We equip $\Sn$ with the Spherical metric, which is the infinitessimal metric 
induced by the Euclidean metric on the ambient space $\R^{n+1}$. 
And we denote by $R=R_n$ the subgroup consisting of Euclidean isometries, 
and by $R_+ := G_+\cap R$ the subgroup of orientation preserving rigid rotations. 
Then $R$ is also the stabilizer of the origin $\zerov$.
For $\wv \in \BBn$ define $g_\wv\in G_+$ by
\EQN{gw}
g_\wv(\xv) = \frac{\xv(1-|\wv|^2) + \wv(1+ |\xv|^2 + 2 <\wv,\xv>)}
{1+ |\wv|^2|\xv|^2+2<\wv,\xv>}, 
\ENDEQN
where $<\cdot,\cdot>$ denotes the Euclidean inner product.
Then $g_\wv$ preserves the line segment $[-\wv/|\wv|, \wv/|\wv|]$, fixes 
$\pm \wv\|\wv|$ and $g_\wv^{-1} = g_{-\wv}$.
Moreover for $0< r <1$ let $\rv = (r, 0, \ldots , 0) = r\ev_1$ and write $g_r := g_\rv$, 
where in general $\ev_j$ denotes the $j$-th element 
of the standard orthonormal basis for $\R^{n+1}$. 
Any $g_\wv$ can be written $g_\wv = \rho\circ g_r$, 
where $r = |\wv|$ and $\rho = g_\wv\circ g_r^{-1}\in R_+$. 
Moreover any element $g\in G_+$ can be written in a unique way as 
$$
g = g_\wv\circ \rho',
$$
where $\wv = g(\zerov)$, 
and 
$$
\rho' = g_\wv^{-1}\circ g = g_{-\wv}\circ g \in R_+.
$$
So that in fact $G_+ = <R_+, (g_r)_{0<r<1}>$ and $G = <R_+, (g_r)_{0<r<1}, c>$.

We can identify $\widehat\R^n$ with $\Sn$ via stereographic projection of the 
central plane $x_{n+1}=0$ in $\widehat\R^{n+1}$ or equivalently through reflection 
in the sphere $\Sn(\ev_{n+1},\sqrt{2})$.
In the case $n=2$ stereographic projection identifies $\Cbar= \CPone$ with $\Sto$. 
And in the $\C$ coordinate an orientation preserving {\Mobius}-transformation $g$ 
preserving the unit circle can be written 
$g = g_w(\rho z)$, where $|\rho| = 1$ and 
$$
g_w(z) = \frac{z + w}{1 + \overline{w}z}, \qquad \textrm{where } w \in \D.
$$

Following Douady and Earle the group $G$ operates on $\BBn$, $\partial\BBn = \Sn$, 
on the set of probability measures $\PP(\Sn)$ and on the vector space $\FF(\BBn)$ 
of continuous vector fields on $\BBn$. That is 
\begin{align*}
&g\cdot\zv = g(\zv), \qquad \textrm{for}~\zv\in \overline{\BBn},\\
&(g\cdot\mu)(A) = g_*\mu(A) = \mu(g^{-1}(A)), \quad \textrm{for}~\mu\in\PP(\Sn)
\textrm{ and } A\subset\Sn\ \textrm{ a Borel subset},\\
&(g\cdot \vv)(g(\zv)) = g_*(\vv)(g(\zv)) = D_\zv g(\vv(\zv)), \quad \textrm{for}~\vv\in\FF(\BBn) 
\textrm{ and } \zv\in \BBn.
\end{align*}
Here $D_\zv g$ denotes the differential of $g$ at $\zv$. 
The group $G\times G$ operates on the spaces $\End(\BBn), \CC(\BBn)$ and $\End(\Sn), \CC(\Sn)$ of endomorphisms and continuous endomorphisms
of $\BBn$ and $\Sn$ respectively by 
$$
(g,h)\phi := g\circ \phi\circ h^{-1}.
$$
If $G$ operates on the spaces $\X$ and $\Y$ then a map {\mapfromto T \X \Y} is called $G$ equivariant or conformally natural if 
$$
\forall\;g\in G,\quad \forall\;x\in \X\quad : T(g\cdot x) = g\cdot T(x).
$$
And if $G\times G$ operates on both $\X$ and $\Y$ then conformal naturality  
of $T$ is taken to mean $G\times G$-equivariance.

Douady and Earle introduced the idea of Conformal Barycenter for probability measures on 
$\Sen\subset\C$ and more generally on $\Sn\subset\R^n$. 
They used the conformal barycenter to define conformally natural extensions 
of self-homeomorphisms of $\Sn$. 
We shall in this note study the application of the Douady-Earle extension operator 
to a much wider class of maps than just homeomorphisms.
More precisely let $\eta_\zerov$ denote the normalized 
standard Euclidean or Lebesgue probability measure on $\Sn$. 
And let {\mapfromto f \Sn \Sn} be a measureable endomorphism, 
for which the push-forward $f_*(\eta_\zerov)$ of $\eta_\zerov$ by $f$ 
is absolutely continuous with respect to $\eta_\zerov$.
We shall show that the Douady-Earle extension operator also yields 
a conformally natural extensions of maps such as $f$, 
to non-constant self maps also denoted by $f$, 
{\mapfromto f {\overline{B_{n+1}}} {\overline{B_{n+1}}}}, 
which are real-analytic in the interior and continuous, 
whenever the original map $f$ is continuous. 
In particular we obtain extensions of rational and entire transcendental maps 
of $\Sto\isom\CPone$ to the hyperbolic three-space $\BBtre$. 
And of course by reflection and renewed stereographic projection to 
an endomorphism of $\Stre$.

The motivation for this note comes from talks by Bill Thurston 
at a workshop in Roskilde 2010, where he asked: 
"What is the three-manifold of a rational map?" or 
"How can we in a natural way extend rational maps to $\BBtre$?". 
The answer we propose to the second question is: Use the Douady-Earle extension.

In order to be self-contained we shall start by reviewing the Douady-Earle construction of 
the conformal barycenter and the Douady-Earle extension in any dimension. 

\section{Conformal Barycenters}
\subsection{Harmonic measure}

Denote be $\eta_0$ the normalized Euclidean Lebesgue measure on $\Sn$, 
$$
\eta_0(A) = \frac{1}{\Vol(\Sn)} \int\ldots\int_A \d L, \qquad \Vol(\Sn) =  \int\ldots\int_A \d L,
$$
where $L$ denotes Lebesgue measure. 
We shall henceforth also write 
$$
\eta_0(A) = \int_A \deta_0.
$$
Then $\eta_0$ is invariant under $R$, i.e. $g_*(\eta_0)=\eta_0$ for every element $g\in R$. 

For $\wv \in \BBn$ the harmonic measure with center $\wv$ is the measure $\eta_\wv = (g_\wv)_*(\eta_0)$. 
Note that by the above $\eta_\wv = g_*(\eta_0)$ for any $g \in G$ with $g(0) = \wv$.

Also note that since each $g \in G$ is conformal 
$$
|\Jac_g(\zv)| = {||\Jac_g(\zv)||}^n,
$$
where $||\cdot||$ denotes the operator norm and $|\cdot|$ denotes determinant.
In the $2$-dimensional and thus $1$-complex dimensional case one computes for $g_w$ and $|z|=1$:
$$
|g_w'(z)| = \frac{1-|w|^2}{{|z+w|}^2}.
$$
Thus in real dimension $n$ we obtain for $\zv\in\Sn$ and $\wv\in \B_{n+1}$:
$$
|\Jac_{g_\wv}(\zv)| = {\left(\frac{1-|\wv|^2}{{|\zv+\wv|}^2}\right)}^n
$$
and hence by the change of variables formula
$$
\eta_\wv(A) = \int_A 1 \deta_\wv = \int_A 1 \d (g_\wv)_*\eta_\zerov  
= \int_{g_{-\wv}(A)} 1 \deta_\zerov = \int_A {\left( \frac{1-|\wv|^2}{{|\zv-\wv|}^2}\right)}^n \deta_\zerov.
$$

\subsection*{The Conformal Barycenter of a measure}
Let us define a probability measure to be \emph{admissible}, 
if it has no atoms of mass greater than or equal to $1/2$. 
And let $\PPp(\Sn)$ denote the space of admissible probability measures.
To each admissible probability measure $\mu\in\PPp(\Sn)$ 
we shall assign a point $B(\mu)\in \BBn$ 
so that the map $\mu \mapsto B(\mu): {\PP'(\Sn)} \to \BBn$ 
is conformally natural and normalized by
\EQN{centered}
B(\mu) = \zerov \quad \Leftrightarrow \quad \int_{\Sn} \zetav \dmu(\zetav) = \zerov
\ENDEQN

\REFPROP{Vdef}
The mapping {\mapfromto V {\PP(\Sn)} {\FF(\BBn)}}, 
which to a probality measure $\mu\in\PPp(\Sn)$ assigns the vector field
\EQN{Vfield}
V_\mu(\wv) = \frac{1- |\wv|^2}{2} \int_\Sn g_{-\wv}(\zetav) \dmu(\zetav), \qquad \wv\in\BBn
\ENDEQN
is the unique conformally natural such map satisfying the normalizing 
condition
\EQN{Vfieldatorigin}
V_\mu(\zerov) = \frac{1}{2}\int_\Sn \zetav \dmu(\zetav).
\ENDEQN
\ENDPROP
The normalizing factor $\frac{1}{2}$ is inessential, but kept here in order to make 
$V_\mu$ asymptotically a hyperbolic unit vector field at $\infty$, 
when $\mu$ has no atoms.
\PROOF{}
Equivariance or conformal invariance is equivalent to
$$
\forall\, g\in G, \forall\,\wv\in\BBn : \qquad V_{g_*\mu}(g(\wv)) = 
(g\cdot V_\mu)(g(\wv)) = 
D_\wv g(V_\mu(\wv))
$$
Thus the normalizing condition \eqref{Vfieldatorigin} is invariant under the subgroup $R$ 
stabilizing the origin, because such maps are linear. 
And for $g=g_{-\wv} = g_\wv^{-1} $ with $g_{-\wv}(\wv) = \zerov$, 
the above formula implies
$$
\forall\, \wv\in\BBn : \qquad V_\mu(\wv) = D_{\zerov}g_{\wv}V_{(g_{-\wv})_*\mu}(\zerov)
$$
Thus the mapping $\mu \mapsto V_\mu$ is conformally natural if and only if
\begin{align}
V_\mu(\wv) &= \frac{1- |\wv|^2}{2} \int_\Sn \zetav \d ((g_{-\wv})_*\mu)(\zetav)\notag\\
&= \frac{1- |\wv|^2}{2} \int_\Sn g_{-\wv}(\zetav) \dmu(\zetav)\notag\\
&= \frac{1- |\wv|^2}{2} \int_\Sn \frac{\zetav(1-|\wv|^2) - 2\wv(1-<\zetav,\wv>)}
{1+|\wv|^2-2<\zetav,\wv>} \dmu(\zetav).\label{Vformula}
\end{align}
\ENDPROOF

Next we want to prove that 
\REFPROP{uniquezero}
For each admidsisble measure $\mu\in \PPp(\Sn)$ 
the vector field $V_\mu$ has a unique zero in $\BBn$.
\ENDPROP
For the proof we shall use a few elementary lemmas, 
which are generalizations to dimension $3$ and higher of the corresponding statements for the complex plane, as can be found in 
\cite[Sections 2 and 11]{D-E}.
 
\REFLEM{stableiquilibrium}
For any admisble probability measure $\mu\in \PPp(\Sn)$ 
any zero $\vv\in\BBn$ of the vector field $V = V_\mu$ 
is an isolated stable equilibrium.
\ENDLEM
\PROOF{}
By conformal naturallity it suffices to consider the case $\vv = \zerov$. 
Expanding the above formula \eqref{Vformula} for $V_\mu(\wv)$ to first order in $\wv$ we obtain:
\begin{align*}
V_\mu(\wv) &= \frac{1}{2}\int_\Sn 
\zetav - 2(\wv-\zetav<\wv,\zetav>) \dmu(\zetav) + \ooo(|\wv|)\\
&= V_\mu(\zerov) - \int_\Sn (\wv-\zetav<\wv,\zetav>) \dmu(\zetav) + \ooo(|\wv|)\\
&= - \int_\Sn (\wv-\zetav<\wv,\zetav>) \dmu(\zetav) + \ooo(|\wv|)\\
\end{align*}
since $V_\mu(\zerov) = \zerov$. Hence the Jacobian of $V$ at $\vv=\zerov$ is given by
\begin{equation}\label{diffVmu} 
\Jac_V(\zerov)(\epsv) = - \int_\Sn (\epsv -\zetav<\epsv,\zetav>)\dmu(\zetav)
\end{equation}
and thus $\Jac_V(\zerov)$ is non singular. In fact $\vv$ is a sink since 
\EQN{Vmuquadformnegative}
<\epsv, \Jac_V(\zerov)(\epsv)> 
= - \int_\Sn (<\epsv,\epsv>-<\zetav,\epsv><\epsv,\zetav>) \dmu(\zetav)
< 0.
\ENDEQN
\ENDPROOF
Douady and Earle showed that if $\mu(D_{\Sn}(\ev_1,\pi/4)) \geq \frac{2}{3}$ then 
\EQN{DEsimplebound}
<V_\mu(0),\ev_1> > 0,
\ENDEQN
where $D_{\Sn}(\ev_1,\delta)$ denotes the closed ball in $\Sn$ of center $\ev_1$ 
and spherical radius $\delta$.
This is sufficient to prove \propref{uniquezero}, 
if $\mu$ has no atoms of mass $\frac 1 3$ or higher. 
To prove the Propositon also, when no atom has mass $\frac 1 2$ or higher, 
we need the following slight refinement:
\REFLEM{directionatorigin}
Let $\delta\in]0,\sqrt{2}[$ and suppose 
$\mu(D_{\Sn}(\ev_1,\delta)) \geq (1+\frac{\delta^2}{2})/2$. 
Then 
$$<V_\mu(0),\ev_1> > 0.
$$
\ENDLEM
\PROOF{}
\begin{align*}
<V_\mu(0),\ev_1> &= \int_{D_{\Sn}(\ev_1,\delta)} <\zetav, \ev_1> \dmu(\zetav) +  
\int_{\Sn\Sm D_{\Sn}(\ev_1,\delta)} <\zetav, \ev_1> \dmu(\zetav)\\
&\geq(1-\frac{\delta^2}{2})(1+\frac{\delta^2}{2})/2 - 1\cdot (1-\frac{\delta^2}{2})/2 = 
\frac{\delta^2}{4}(1-\frac{\delta^2}{2}) > 0.
\end{align*}
\ENDPROOF
\REFLEM{outerinwards}
Suppose that $\mu$ is admissible. 
Then there exists $r\in\;]0,1[$ such that $V_\mu(\wv)$ points inwards 
at any point $\wv\in\BBn$ with 
$r \leq |\wv| < 1$, i.e.~$<V_\mu(\wv),\wv> < 0$. 
\ENDLEM
\PROOF{}
Choose $\delta \in\; ]0,\sqrt{2}[$ such that for any $\zetav\in\Sn$ : 
$\mu(\{\zetav\}) < (1-\frac{\delta^2}{2})/2$. 
Then there exists $\epsilon\in\; ]0,\pi[$ such that for  $\zetav\in\Sn$ : 
$\mu(D_{\Sn}(\zetav,\eps)) \leq (1-\frac{\delta^2}{2})/2$. 
Choose $r\in\;]0,1[$ such that $\forall\;\wv\in\BBn$ with $r\leq |\wv| < 1$:
$$
\eta_\wv(\Sn\Sm D_{\Sn}(\frac \wv {|\wv|},\epsilon)) \leq \eta_\zerov(D_{\Sn}(\ev_1,\delta))
$$
Then it follows from \lemref{directionatorigin} that $V_\mu(\wv)$ 
points into the sphere $ S = |\wv|\Sn$: 
Let $g\in G_+$ be any {\Mobius} transformation mapping $\wv$ to $\zerov$ 
and $-\wv/|\wv|$ to $\ev_1$. 
Then $g(S)$ is a sphere through $\zerov$ 
and with $\ev_1$ as an inwards pointing normal vector at $\zerov$.
Moreover let $\nu = g_*\mu$ then by conformal naturality $g_*(V_\mu(\wv)) = V_\nu(0)$ 
and $\nu$ satisfies the hypotheses of \lemref{directionatorigin}. 
\ENDPROOF

\PROOF{of \propref{uniquezero}}
Let $\mu\in\PP(\Sn)$ be any admissible measure, i.e.~with no atom of mass $1/2$ or higher. 
In \lemref{stableiquilibrium} we have shown that any zero of the vector field $V_\mu$ 
is an isolated stable equilibrium, i.e.~the vector field points inwards on small spheres around the zero. 
Moreover by \lemref{outerinwards} the vector field $V_\mu$ 
is pointing inwards near the boundary $\Sn$ of $\BBn$. 
Hence by the \Poincare-Hoppf theorem \cite[see also Lemma 3, p 36]{Milnor} $V_\mu$ has a unique zero $B(\mu)\in\BBn$.
\ENDPROOF

\REFDEF{conformalbarycenter}
Define a conformally natural mapping {\mapfromto B {\PPp(\Sn)} {\BBn}} by 
setting $B(\mu)$ equal to the unique zero $\wv\in\BBn$ of the vector field $V_\mu$. 
Then $B$ satisfies \eqref{centered}.
\ENDDEF

\section{Extending continuous endomorphisms of $\Sn$.}
Let $\EE(\Sn)$ denote the space of endomorphisms {\mapfromto \phi \Sn \Sn} 
such that $\phi_*\eta_\zerov$ has no atoms. 
For such mappings the measures $\phi_*\eta_\zv$ has no atoms neither for any $\zv\in\BBn$.
And let $\End(\BBnbar)$ denote the space of endomorphisms of $\BBnbar$, 
whose restrictions to $\BBn$ are endomorphisms of $\BBn$.

The Douady-Earle extension operator $E$, in the following denoted the D-E extension,  
which Douady and Earle studied for homeomorphisms is the map
{\mapfromto E {\EE(\Sn)} {\End(\BBnbar)}} defined as follows:
For $\phi\in\EE(\Sn)$ the mapping 
{\mapfromto {E(\phi) = \Phi} {\overline{\BBn}} {\overline{\BBn}}} 
is given by the formulas
\begin{equation}
\Phi(\zv) = 
\begin{cases}
\phi(\zv),\quad &\zv\in\Sn,\\
B({(\phi\circ g_\zv)}_*(\eta_0)) = B(\phi_*(\eta_\zv)), \quad &\zv\in\BBn
\end{cases}
\end{equation}
Clearly the mapping $\phi\mapsto E(\phi) = \Phi$ is conformally natural, i.e.~
for all $g,h\in G$:
$$
E(g\circ\phi\ h) = g \circ E(\phi) \circ h.
$$
Moreover for any conformal automorphism $g\in G$ we have $E(g_{|\Sn}) = g$, 
by conformal naturality of $E$ and the fact that $B(\eta_\zerov)=\zerov$. 
We can also formulate this as saying that the D-E extension operator extends 
the {\Poincare} extension operator.
For $n=1$ at least we have a much stronger property:
For inner functions, \cite[Def.~17.14]{Rudin} that is for holomorphic selfmaps {\mapfromto f \D \D} 
of the unit disc $\D\subset\C$ with boundary values in $\Sen$ a.e.,  
the D-E extension simply recovers $f$ from its boundary values. 
More precisely it is well known that for bounded holomorphic functions 
(see \cite[Th.~11.21]{Rudin}) the radial limit
$$
f^\#(\zeta) = \lim_{r\nearrow 1}f(r\zeta)
$$
exists for a.e.~$\zeta\in\Sen$ and satisfies the Cauchy formula:
\EQN{Cauchyformula}
\forall\;z\in\D:\qquad f(z) = \frac{1}{2\pi i}\int_\Sen \frac{f^\#(\zeta)}{\zeta-z} \dzeta.
\ENDEQN
Moreover the space of such functions $f^\#$ 
is the space of bounded measureable functions, whose negative Fourier coefficients 
are all equal to zero. 
For inner functions where $|f^\#(\zeta)|=1$ a.e.~the measure 
$f_*^\#\eta_0$ is absolutely continuous with respect to $\eta_0$ 
(see \cite[Th.~17.13]{Rudin}) and hence 
$f_*^\#\eta_z$ is absolutely continuous with respect to $\eta_0$ for any $z\in\D$, 
so that $f^\#\in\EE(\Sen)$.
\REFPROP{Hinnerextensions}
If {\mapfromto f \D \D} is an inner function then
$$
E(f^\#)(z) = f(z),\qquad \forall z\in\D.
$$
\ENDPROP
\PROOF{}
Let $f$ be an arbitrary inner function and let $z\in\D$ be arbitrary. 
We need to show that $f(z)= E(f^\#(z)$. By conformal naturality 
we can suppose $z = f(z) = 0$ as we may precompose by $g_z$ 
and postcompose by $g_{-f(z)}$.
That is it suffices to prove that $B(f_*\eta_0) = 0$ 
for any inner function $f$ with $f(0)=0$.
We compute
$$
2V_{f_*^\#\eta_0}(0) = \int_\Sen \zeta \d(f_*^\#\eta_0)(\zeta)
=\int_\Sen f^\#(\zeta) \deta_0(\zeta)
=\frac{1}{2\pi i}\int_\Sen \frac{f^\#(\zeta)}{\zeta} \dzeta = f(0) = 0.
$$
\ENDPROOF

\REFLEM{continuityofextension}
Let $\phi\in\EE(\Sn)$ and let $\Phi = E(\phi)$. If $\phi$ is continuous at some point
$\zetav_0\in\Sn$ then so is $\Phi$. In particular if $\phi$ is continuous then $\Phi$ 
is continuous on $\Sn$.
\ENDLEM
We shall see in the next lemma that $\Phi$ is real-analytic in $\BBn$, 
so that in particular $\Phi$ is continuous whenever $\phi$ is continuous. 
\PROOF{}
Recall that Euclidean balls $\BBn(\zetav, r)$ and spheres $\Sn(\zetav,r)$ are mapped 
to such balls and spheres (possibly half spaces and hyperplanes or complements of a closed 
ball union $\infty$) under any conformal automorphism $g\in G$.

Thus given a spherical ball $B_\Sn(\zetav,\delta)\subset \Sn$, $0<\delta < \pi$ and $\wv\in\BBn$ 
there are two alternative ways of describing the size of 
$B_\Sn(\zetav,\delta)$ viewed from $\wv$. 
Either we can use the visual Poincar{\'e} radius from $\wv$, 
i.e.~the spherical radius of the ball $g_{-\wv}(B_\Sn(\zetav,\delta))$ in $\Sn$ 
or we can use the $\wv$-harmonic measure 
$\eta_\wv(B_\Sn(\zetav,\delta)) = \eta_\zerov(g_{-\wv}(B_\Sn(\zetav,\delta)))$.

Given $B_\Sn(\zetav,\delta)$ we denote by $W(B_\Sn(\zetav,\delta))$ 
the set $B_\Sn(\zetav,\delta)$ itself union the open subset of points $\wv\in \BBn$ for 
which the visual Poincar{\'e} radius from $\wv$ exceeds $\pi/4$. 
Similarly we denote by $U(B_\Sn(\zetav,\delta))$ the set $B_\Sn(\zetav,\delta)$ itself
union the open subset of points $\wv\in \BBn$ for 
which $\eta_\wv(B_\Sn(\zetav,\delta)) > 2/3$. 
Then $U(B_\Sn(\zetav,\delta)) \subset W(B_\Sn(\zetav,\delta))$ 
and both sets are neighborhoods of $\zetav$ in $\BBnbar = \BBn\cup\Sn$. 
In fact for $\zetav =\ev_1$ and $\delta = \pi/4$ the set $W(B_\Sn(\zetav,\delta))$ 
equals the intersection of $\BBnbar$ with the open ball $\BBn(\sqrt{2}\ev_1, 1)$ 
and for $\delta=2\pi/3$ the upen set $U(B_\Sn(\zetav,\delta))$ 
is the complement $\BBnbar\Sm\BBnbar(-2\ev_1,\sqrt{3})$. 
Clearly any of the families of sets $U(B_\Sn(\zetav,\delta)), W(B_\Sn(\zetav,\delta))$, 
$0<\delta<\pi$ forms fundamental systems of neighbourhoods of $\zetav$ in $\BBnbar$. 
Suppose $\phi$ is continuous at $\zetav_0$ and let $0<\eps< \pi$ be given. 
Choose $0<\delta<\pi$ such that 
$$
\phi(B_\Sn(\zetav_0,\delta)) \subset B_\Sn(\phi(\zetav_0,\eps)).
$$
Then for any $\wv\in U(B_\Sn(\zetav_0,\delta))$ 
and any $\zv\in \partial W(B_\Sn(\phi(\zetav_0,\eps)))\cap\BBn$ the vector 
$V_{\phi_*\eta_\wv}(\zv)$ points into $W(B_\Sn(\phi(\zetav_0,\eps)))$. 
Hence $\Phi(\wv)$ the unique zero of $V_{\phi_*\eta_\wv}$ belongs to 
$W(B_\Sn(\phi(\zetav_0,\eps)))$. This proves continuity at $\zeta_0$.
\ENDPROOF

\REFLEM{internalrealanalyticity}
Let $\phi\in\EE(\Sn)$ and $E(\phi) = \Phi$ be as above. 
Then $\Phi$ is real-analytic in $\BBn$.
\ENDLEM
\PROOF{}
Towards real-analyticity of $\Phi$ recall that $\Phi(\zv)$ is the unique zero 
of the vector field 
\begin{align*}
V_{\phi_*(\eta_\zv)}(\wv) 
&= \frac{1 - {|\wv|}^2}{2}\int_\Sn \zetav \d(g_{-\wv}\circ\phi)_*\eta_\zv(\zetav)\\
&= \frac{1 - {|\wv|}^2}{2}\int_\Sn  g_{-\wv}(\phi(\zetav))\d(g_{-\zv})_*\eta_\zerov(\zetav)\\
&= \frac{1 - {|\wv|}^2}{2}\int_\Sn  g_{-\wv}(\phi(\zetav))
{\left(\frac{1-{|\zv|}^2}{{|\zv-\zetav|}^2}\right)}^n\deta_\zerov(\zetav).
\end{align*}
Thus $\forall\;\zv\in\BBn$ the value $\wv = \Phi(\zv)$ is the unique point $\wv\in\BBn$ 
such that:
\[
F(\zv,\wv) = \frac{2V_{\phi_*(\eta_\zv)}(\wv)}{1-|\wv|^2} 
=\int_\Sn  g_{-\wv}(\phi(\zetav))
{\left(\frac{1-{|\zv|}^2}{{|\zv-\zetav|}^2}\right)}^n\deta_\zerov(\zetav) = \zerov.
\]
Clearly $F$ is a real-analytical function of $(\zv,\wv)\in \BBn\times\BBn$. 
Thus by the implicit function theorem 
we need only show that for any pair $(\zv,\wv)\in \BBn\times\BBn$ 
with $F(\zv,\wv) = \zerov$ the $\wv$ partial derivatives matrix
$\J_\wv F = \frac{\partial F}{\partial\wv}$ evaluated at $(\zv,\wv)$ is non-singular. 
By conformal naturality we can suppose $\zv = \wv = \Phi(\zv) =\zerov$. 
A straight forward computation analogous to the one leading to \eqref{diffVmu} 
yields that $\J_\wv F$ evaluated at $(\zerov,\zerov)$ and 
applied to the vector $\epsv$ is given by the formula:
$$
\J_\wv F(\epsv) 
= -2 \int_\Sn (\epsv-<\epsv,\phi(\zeta)>\phi(\zeta)) \deta_\zerov(\zeta).
$$
Similarly as for $\Jac_V$ this shows that $J_\wv F$ is non singular at $(\zerov,\zerov)$. 
So $\Phi$ is real-analytic by the implicit function theorem and
\EQN{JacPhi}
\Jac_\Phi(\zerov) = -{(\J_\wv F)}^{-1}\circ \J_\zv F
\ENDEQN
where both partial derivative matrices are evaluated at $(\zerov,\zerov)$ and 
$$
\J_\zv F = \int_\Sn \phi(\zetav)\times\zetav \deta_\zerov(\zetav),
$$
and where $\phi(\zetav)\times\zetav$ is the matrix valued mapping
$$
A_{i,j}(\zetav) = \phi_i(\zetav)\cdot \zeta_j.
$$
\ENDPROOF

\section{Properties of the D-E extension of rational maps}
The question that naturally arises is: For $f$ a rational map on the Riemann sphere. 
What are the geometric and dynamical properties of the D-E extension $E(f)$?
How many of the properties of $f$ are inherited by $E(f)$.
By elementary topology $E(f)$ is a proper map, 
that is the preimage of any compact set is compact. 
And moreover for any point $\wv\in\overline\BBn$ 
the pre image $E(f)^{-1}(\wv)$ is a real analytic set.\\
\\
\noindent Question 1: Is $E(f)$ a discrete map?\\
Question 2: Is $E(f)$ an open map?\\
Question 3: Is $E(f)$ a map of the same degree as $f$? \\
Question 4: Is the Julia set of $E(f)$ (the set of points $\xv$ 
for which the family of iterates does not form an equicontinuous family on any neighbourhood of 
$\xv$) equal to the convex hull of the Julia set for $\wf$?\\

In certain elementary cases at least the immediate answer to the above questions are yes, but not completely satisfactory.

For the following discussion we shall identify $\C$ with the coordinate plane in $\R^3$, 
$\{\xv=(x_1,x_2,x_3)|\; x_3=0 \}$ and write $z=x+iy$ for the point $(x,y,0)$. 
In particular we shall identify the complex unit disk $\D$ with the disc 
$\{\xv\in\R^3|\; |z|^2 = x_1^2+x_2^2 <1, x_3=0\}$ and 
the unit circle $\Sen$ with the circle $\{\xv\in\R^3|\; |z|^2 = 1, x_3=0\}$. 
Then stereographic projection $S$ of $\Cbar$ on to $\Sto$ from the north pole 
$N = \ev_3\in\R^3$ is the map
$$
z \mapsto S(z) = \left(\frac{2z}{1+|z|^2}, \frac{1-|z|^2}{1+|z|^2}\right) = 
\frac{1}{1+|z|^2}(2x,2y, 1-|z|^2).
$$
For {\mapfromto f {\Cbar\;(\C)} \Cbar} a holomorphic map we shall write 
$\wf$ for its conjugate by $S$, i.e.:
$$
\wf(S(z)) = S(f(z)).
$$
In the following we shall discuss finite Blaschke products 
$$
f(z) = \sigma\prod_{j=1}^d\frac{z+a_j}{1+\overline{a}_jz},
\qquad|\sigma|=1, \quad a_j\in\D
$$
\PROP
For $f$ a finite Blaschke product the D-E extension $E(\wf)$ maps $\D$ 
onto $\D$, preserves the upper and lower hemispheres, $\Sto_+, \Sto_-$ and 
further more on $\D$ 
we have $\partial E(\wf)/\partial x_3 = g(z) \ev_3$ for some positive 
real analytical function {\mapfromto g \D \Rplus}.

If moreover $f(z) = z^d$ (i.e.~$a_j=0$ for all j), then $E(\wf)(z) = z^d\cdot h(|z|^2)$ for 
some real analytical function $h$ with $h(r) \to 1$ as $r\nearrow 1$. 
\ENDPROP
\PROOF
The reflection $c(x_1,x_2,x_3) = (x_1,x_2,-x_3)$ is the {\Poincare} extension of $\widehat\tau$, 
where $\tau(z) = 1/\overline{z}$ denotes the reflection in $\Sen$. 
Then $c\circ\wf\circ c = \wf$.
Write $\Phi = E(\wf)$ for the D-E extension of $\wf$.
Then by conformal naturality of the D-E extension 
\EQN{reflsym}
c\;\circ \Phi = \Phi\circ c.
\ENDEQN 
Hence $\Phi(\D)\subseteq \D$, since {\mapfromto {c_|} \D \D} is the identity. 
Moreover if $\Phi(\D) \not= \D$, then a simple homotopy argument would imply that 
the restriction $\wf_| = \Phi_|$ to $\Sen$ is homotopic to a constant map. 
Thus $\Phi(\D) = \D$. 

To prove that $\Phi$ preserves the upper and hence the lower hemisphere it 
suffices to prove that for any $\xv\in\Sto_+$ and any $w\in\D$: 
$\ev_3\cdot V_{\wf_*\eta_\xv}(w) > 0$. 
Furthermore by conformal naturality it suffices to consider the case 
$\xv = t\ev_3$ with $0<t<1$ and $w=0$. 
Before we start computing let us note that since $c_*\eta_\zerov = \eta_\zerov$ 
we have for any measurable function {\mapfromto \phi \Sn {\R\;(\C)}}
$$
\int_\Sn \phi(\zetav) \deta_\zerov = 
\int_\Snplus (\phi(\zetav) + \phi(c(\zetav)))\deta_\zerov(\zetav) = 
\int_\Snminus (\phi(\zetav) + \phi(c(\zetav)))\deta_\zerov(\zetav).
$$
Applying this to $V_{\wf_*\eta_\xv}(0)$ we obtain
\begin{align*}
\ev_3\cdot V_{\wf_*\eta_\xv}(0) &= \int_\Sto \ev_3\cdot \wf(\zetav) 
{\left(\frac{1-|\xv|^2}{|\xv-\zetav|^2}\right)}^2 \deta_\zerov(\zetav)\\
&= \int_{\Sto_+}\left(\wf_3(\zetav){\left(\frac{1-|\xv|^2}{|\xv-\zetav|^2}\right)}^2+
\wf_3(c(\zetav)){\left(\frac{1-|\xv|^2}{|\xv-c(\zetav)|^2}\right)}^2\right)\deta_\zerov(\zetav)\\
&= \int_{\Sto_+}\wf_3(\zetav)\left({\left(\frac{1-|\xv|^2}{|\xv-\zetav|^2}\right)}^2-
{\left(\frac{1-|\xv|^2}{|\xv-c(\zetav)|^2}\right)}^2\right)\deta_\zerov(\zetav) > 0,
\end{align*}
since $|\xv-\zetav| < |\xv-c(\zetav)|$ and $\wf_3$ is positive on $\Sto_+$.

To compute the partial derivative vector $\partial\Phi/\partial x_3(z)$ for $z\in\D$ 
we equate the Jacobians of the two sides of 
\eqref{reflsym} and obtain
$$
\frac{\partial\Phi_1}{\partial x_3}(z) = \frac{\partial\Phi_2}{\partial x_3}(z) = 
\frac{\partial\Phi_3}{\partial x_1}(z) = \frac{\partial\Phi_3}{\partial x_2}(z) = 0.
$$
Thus $\partial\Phi/\partial x_3(z) = \partial\Phi_3/\partial x_3(z) \ev_3 = g(z)\ev_3$. 
To complete the first set of statements we just need to show that $g$ is a positive function. 
By conformal naturality is suffices to consider the case $\Phi(0) = 0$. 
Furthermore in order to simply notation let us for {\mapfromto \phi \Sn {\R\;(\C)}} 
a measureable function write 
$$
M(\phi) := \int_\Sn \phi(\zetav) \deta_\zerov(\zetav).
$$
Then an elementary calculation, using \eqref{JacPhi} yields
$$
g(0) = \frac{\partial\Phi_3}{\partial x_3}(0) = 
\frac{M(\wf_3\cdot \zeta_3)}{2M(1-\wf_3^2)} > 0,
$$
since $\wf_3\cdot \zeta_3, (1-\wf_3^2) \geq 0$ with equality for $\zeta_3=0$ 
only in the first and for $f(z)$ equal to $0$ or $\infty$ in the second. 

In the special case $f(z) = z^d$ we have $f(\e^{i\theta}z) = \e^{id\theta}f(z)$. 
and thus by conformal naturality $\Phi(\e^{i\theta}z) = \e^{id\theta}\Phi(z)$. 
Hence $\Phi(0) = 0$ and for $z\not= 0$:
$$
\Phi(z) = \Phi(|z|)\frac{z^d}{|z|^d}.
$$ 
Moreover $f$ commutes with complex conjugation, which translates to 
$\wf$ commutes with the reflection $(x_1,x_2,x_3) \mapsto  (x_1,-x_2,x_3)$. 
As above this implies that $\Phi(]-1,1[)\subseteq ]-1,1[$, 
so that $\Phi$ is a real analytic real function on the reals.
Expanding the real-analytic function in a power series in $z,\zbar$ on a neighbourhood of 
$0$ and noting that $z\zbar = |z|^2$ we obtain:
$$
\Phi(z) = \frac{z^d}{|z|^d} \sum_{m=0}^{\infty} b_m|z|^m = 
\sum_{n,k=0}^\infty a_{n,k}z^n\zbar^k.
$$
By the uniqueness theorem for power series this implies that
$a_{n,k}=0$ for $n-k\not= d$. Thus we are left with
$$
\Phi(z) = \sum_{k=0}^\infty a_{k+d,k}z^{k+d}\zbar^k = z^d\sum_{k=0}^\infty a_{k+d,k}|z|^{2k}.
$$
\ENDPROOF
Write $M_t(z) = tz$ for $0<t$ so that $M_t$ is a homothety. 
Then the conformal automorphism $h_t = E(\widehat M_t) = g_\wv$ with 
$\wv =\wv(t) = \frac{t-1}{t+1}\ev_3$ maps $\D$ conformally onto the geodesic disk $\D_t$ 
in $\BBtre$ with boundary the circle $\widehat M_t(\Sen)$.
\COR
For $f(z) = z^d$ the D-E extension $E(\wf)$ maps $\D_t$ onto $\D_{t^d}$ by 
a degree $d$ ramified covering and the interval $[\zerov,\ev_3[$ onto itself by an increasing diffeomorphism.
\ENDCOR

\CONJ
For all finite Blaschke products $f$ we have $f=E(\wf)$ on $\D$.
\ENDCONJ
By conformal naturallity of the D-E extension this conjecture is 
equivalent to the seemingly simpler conjecture:
\CONJ
For all finite Blaschke products $f$ with $f(0) = 0$ we have $E(\wf)(\zerov) = \zerov$.
\ENDCONJ
If true the following stronger conjecture would yield almost complete topological understanding 
of $E(\wf)$ for any finite Blaschke product $f$:
\CONJ
For all finite Blaschke products $f$ with $f(0) = 0$ the D-E extension $E(\wf)$ 
maps the interval $[\zerov,\ev_3[$ diffeomorphically and increasingly onto itself.
\ENDCONJ
Clearly the last conjecture implies the two previous ones. 
Moreover by conformal naturality of the D-E extension, it would imply 
that for each $z\in\D$ the unique hyperbolic geodesic through $z$ 
and orthogonal to $\D$ would be mapped diffeomorphically onto the unique such geodesic through 
$f(z)$. And thus the dynamics of $E(\wf)$ would be conjugate to a skew product on 
$\D\times[-1,1]$. Which would be completely analogous to the case of Fuchsian groups.


Address:

Carsten Lunde Petersen, NSM,
Roskilde University,
Universitetsvej 1,
DK-4000 Roskilde,
Denmark.
e-mail: lunde@ruc.dk

\end{document}